\documentclass[12pt]{article}
\usepackage{amssymb,amsmath}
\usepackage{cases}
\usepackage{amsfonts}
\usepackage{color,xcolor}
\usepackage[left=2.0cm,right=2.0cm,top=2.0cm,bottom=2.0cm]{geometry}
\usepackage[colorlinks,citecolor=blue,urlcolor=blue]{hyperref}

\newtheorem{theorem}{Theorem}[section]

\newtheorem{lemma}{Lemma}[section]
\newtheorem{proposition}{Proposition}[section]

\newtheorem{definition}{Definition}[section]
\newtheorem{remark}{Remark}[section]

\newcommand{\bal}{\begin{align}}
\newcommand{\bbal}{\begin{align*}}
\newcommand{\beq}{\begin{equation}}
\newcommand{\eeq}{\end{equation}}
\newcommand{\bca}{\begin{cases}}
\newcommand{\eca}{\end{cases}}

\newcommand{\pa}{\partial}

\newcommand{\dd}{\mathrm{d}}

\newcommand{\R}{\mathbb{R}}

\newcommand{\Z}{\mathbb{Z}}

\newcommand{\bi}{\Big}

\begin{document}
\title{Non-uniform continuity of the generalized Camassa-Holm equation in Besov spaces}

\author{Jinlu Li$^{1}$, Xing Wu$^{2}$\footnote{E-mail: lijinlu@gnnu.edu.cn; ny2008wx@163.com(Corresponding author); mathzwp2010@163.com; 1612921989@qq.com}, Weipeng Zhu$^{3}$ and Jiayu Guo$^{1}$\\
\small $^1$\it School of Mathematics and Computer Sciences, Gannan Normal University, \\
\small Ganzhou, Jiangxi, 341000, China\\
\small $^2$\it College of Information and Management Science,
Henan Agricultural University,\\
\small Zhengzhou, Henan, 450002, China\\
\small $^3$\it School of Mathematics and Information Science, Guangzhou University,\\
\small  Guangzhou, Guangdong, 510006, China}

\date{}

\maketitle\noindent{\hrulefill}

{\bf Abstract:} In this paper, we consider the Cauchy problem for the generalized Camassa-Holm equation proposed by Hakkaev and Kirchev (2005) \cite{Hakkaev 2005}. We prove that the solution map of the generalized  Camassa-Holm equation  is not uniformly continuous
on the initial data in Besov spaces. Our result include the present work (2020) \cite{Li 2020,Li 2020-1} on Camassa-Holm equation with $Q=1$ and extends the previous non-uniform continuity in Sobolev spaces (2015) \cite{Mi 2015} to Besov spaces. In addition, the non-uniform continuity in critical space $B_{2, 1}^{\frac{3}{2}}(\mathbb{R})$ is the first to be considered in our paper.

{\bf Keywords:} Generalized Camassa-Holm equation, Non-uniform continuous dependence, Besov spaces

{\bf MSC (2010):} 35B30; 35G25; 35Q53
\vskip0mm\noindent{\hrulefill}

\section{Introduction}

In this paper, we are concerned with the Cauchy problem for the generalized Camassa-Holm equation introduced by Hakkaev and Kirchev \cite{Hakkaev 2005}
\begin{eqnarray}\label{eq1}
        \left\{\begin{array}{ll}
         u_t-u_{xxt}+\frac{(Q+2)(Q+1)}{2}u^Qu_x=\big(\frac{Q}{2}u^{Q-1}u_x^2+u^Qu_{xx}\big)_x, ~~t>0, ~x\in \mathbb{R},\\
          u(0, x)=u_0, ~~x\in \mathbb{R},\end{array}\right.
        \end{eqnarray}
where $Q\geqslant 1$ is a positive integer and $u(t, x)$ stands for the fluid velocity at time $t\geqslant 0$ in the spatial direction.

Hakkafv and Kirchev \cite{Hakkaev 2005} proved that the Cauchy problem for (\ref{eq1}) is locally well-posed in Sobolev space $H^s(\mathbb{R})$ with $s>\frac{3}{2}$. The local well-posed in the sense of  Hadamard in Besov space $B_{p, r}^s(\mathbb{R})$, $s>\max\{1+\frac{1}{p}, \frac{3}{2}\}$, $1\leq p, r\leq\infty$ or $B_{2, 1}^\frac{3}{2}(\mathbb{R})$ were established in \cite{Mi 2013, Yan 2014, 1Yan 2014}, respectively. Mi and Mu \cite{Mi 2015} showed  that the Cauchy problem for (\ref{eq1}) is locally well-posed in Sobolev space $H^s$ with $s>\frac{3}{2}$ for both the periodic and the non-periodic case. In addition, they proved that the solution map is not uniformly continuous.

 When $Q=1$, Eq. (\ref{eq1}) is reduced to the famous Camassa-Holm (CH) equation
 \begin{eqnarray}\label{eq2}
          u_t-u_{xxt}+3uu_x=2u_xu_{xx}+uu_{xxx},
                  \end{eqnarray}
which was first found by Fuchssteiner and Fokas \cite{Fokas 1981} and later derived as a water wave model by Camassa and Holm \cite{Camassa 1993}. The CH equation is completely integrable \cite{Camassa 1993, Constantin 2001},  has bi-Hamiltonian structure \cite{Fokas 1981, Constantin 1997}, possesses an infinity of conservation laws and admits exact peakon solutions of the form $ce^{-|x-ct|}$ , describing an essential feature of the travelling waves of largest amplitude \cite{Camassa 1993, Constantin 2007, Constantin 2011}.

 The local well-posedness and
ill-posedness for the Cauchy problem of the CH equation (\ref{eq2}) in Sobolev spaces and Besov spaces were studied in the
series of papers  \cite{Constantin 1998, 1Constantin 1998, Danchin 2001, Guo 2019, Li 2016}. After the non-uniform continuity for some dispersive equations was studied by Kenig et al. \cite{Kenig 2001}, non-uniform dependence of the CH equation has been studied by several authors.
The first non-uniform dependence result for CH equation was established by Himonas and Misio{\l}ek \cite{Himonas 2005} in $H^s$ with $s\geqslant 2$ on the circle using explicitly constructed travelling wave
solutions and this was then  sharpened to $s>\frac{3}{2}$ in \cite{Himonas 2009} on the line and \cite{Himonas 2010} on the circle. The above mentioned works utilize the method of approximate solutions in conjunction with delicate commutator and multiplier estimates.
Danchin\cite{Danchin 2001, Danchin 2003} showed the local existence and uniqueness of strong solutions to (\ref{eq2}) with the initial data in $B_{p, r}^s$ for  $s>\max\{1+\frac{1}{p}, \frac{3}{2}\}$, $1\leq p\leq \infty$, $1\leq r<\infty$. The continuous properties of the solutions on the initial data have been supplemented by Li and Yin in \cite{Li 2016}. Very recently, Li, Yu and Zhu \cite{Li 2020} have sharpened the results in \cite{Li 2016}  by showing that the solution map is not uniformly continuous, which depends deeply on the estimates of the transfer equation in Besov space and the constructed high-low frequency smooth initial data.

For studying the non-uniform continuity of the generalized Camassa-Holm equation, it is
more convenient to express (\ref{eq1}) in the following equivalent nonlocal form

\begin{equation}\label{eq3}
\begin{cases}
u_t+u^Q\pa_xu=-\pa_x(1-\pa^2_x)^{-1}[\frac{Q^2+3Q}{2(Q+1)}u^{Q+1}+\frac Q2u^{Q-1}(\pa_xu)^2]:=\mathbf{P}(u), \; &(t,x)\in \R^+\times\R,\\
u(0,x)=u_0,\; &x\in \R.
\end{cases}
\end{equation}

Up to the present, there is no result for the non-uniform continuous dependence of Eq. (\ref{eq3}) in Besov space and it seems more difficult due to the fact that the structure of (\ref{eq3}) when $Q\geq 2$ is much more complicated than that of the CH equation. On the one hand, when bounding the core part $u-u_0+t(u_0)^Q\partial_xu_0$ by $t^2$ in the small time near $t=0$, because of the high nonlinearity in convection term, the method developed for the CH equation in \cite{Li 2020} will only aggravate this difficulty. To overcome this difficulty, we add a new term $t\mathbf{P}(u_0)$, and by using the differential mean value theorem, it is transformed into the estimates of nonlocal term $\mathbf{P}(u)-\mathbf{P}(u_0)$ and convective term $u^Q\partial_xu-(u_0)^Q\partial_xu_0$, which is valid for any bounded set in the above working Besov space. Thus, for any bounded set $u_0$ in working space, the corresponding solution $\mathbf{S}_t(u_0)$ can be approximated by a function of first degree of time $t$ with respect to initial data $u_0$. That is,
$$\mathbf{S}_t(u_0)=u_0-t(u_0)^Q\partial_xu_0+t\mathbf{P}(u_0)+t^2\mathcal{O}(u_0),$$
where $\mathcal{O}(u_0)$ is a bounded quantity of some Besov norm of  $u_0$. Therefore, the non-uniform continuity of $\mathbf{S}_t(u_0)$ is transformed into that of $u_0-t(u_0)^Q\partial_xu_0+t\mathbf{P}(u_0)$, which is the essence in our paper and is different from that in \cite{Li 2020}. On the other hand,  more precise product and multiplier estimates are necessary. In addition, the non-uniform continuity in critical space $B_{2, 1}^{\frac{3}{2}}(\mathbb{R})$ is  considered for the first time.

Our main result is stated as follows.
\begin{theorem}\label{the1.1} Let $s>\max\{1+\frac{1}{p}, \frac{3}{2}\}$, $1\leq p, r\leq \infty$ or $(s,p,r)=(\frac32,2,1)$. The solution map $u_0\rightarrow \mathbf{S}_t(u_0)$ of the initial value problem (\ref{eq3}) is not
uniformly continuous from any bounded subset of  $B_{p, r}^s(\mathbb{R})$ into $\mathcal{C}([0, T];  B_{p, r}^s(\mathbb{R}))$. More precisely,
there exist two sequences $g_n$ and $f_n$ such that
\begin{eqnarray*}
        \|f_n\|_{B_{p, r}^s}\lesssim 1,\qquad  \qquad \qquad \lim_{n\rightarrow \infty} \|g_n\|_{B_{p, r}^s}=0,
        \end{eqnarray*}
but
\begin{eqnarray*}
      \liminf_{n\rightarrow \infty} \|\mathbf{S}_t(f_n+g_n)-\mathbf{S}_t(f_n)\|_{B_{p, r}^s}\gtrsim t, \qquad t\in [0, T_0],
        \end{eqnarray*}
with small positive time $T_0$ for $T_0\leq T$.
\end{theorem}

\begin{remark}\label{rem1}
Our result include the present work \cite{Li 2020,Li 2020-1} on Camassa-Holm equation with $Q=1$.
\end{remark}

\begin{remark}\label{rem2}
Since $B_{2, 2}^s=H^s$, our result covers and extends the previous non-uniform continuity in Sobolev space \cite{Mi 2015} to Besov space.
\end{remark}

{\bf Notations}:  Given a Banach space $X$, we denote the norm of a function on $X$ by $\|\|_{X}$, and \begin{eqnarray*}
\|\cdot\|_{L_T^\infty(X)}=\sup_{0\leq t\leq T}\|\cdot\|_{X}.
\end{eqnarray*} The symbol
$A\lesssim B$ means that there is a uniform positive constant $C$ independent of $A$ and $B$ such that $A\leq CB$.

\section{Littlewood-Paley analysis}

In this section, we will review the definition of Littlewood-Paley decomposition and nonhomogeneous Besov space, and then list some useful properties. For more details, the readers can refer to \cite{Bahouri 2011}.

There exists a couple of smooth functions $(\chi,\varphi)$ valued in $[0,1]$, such that $\chi$ is supported in the ball $\mathcal{B}\triangleq \{\xi\in\mathbb{R}^d:|\xi|\leq \frac 4 3\}$, $\varphi$ is supported in the ring $\mathcal{C}\triangleq \{\xi\in\mathbb{R}^d:\frac 3 4\leq|\xi|\leq \frac 8 3\}$ and $\varphi\equiv 1$ for $\frac{4}{3}\leq |\xi| \leq \frac{3}{2}$. Moreover,
$$\forall\,\ \xi\in\mathbb{R}^d,\,\ \chi(\xi)+{\sum\limits_{j\geq0}\varphi(2^{-j}\xi)}=1,$$
$$\forall\,\ \xi\in\mathbb{R}^d\setminus\{0\},\,\ {\sum\limits_{j\in \Z}\varphi(2^{-j}\xi)}=1,$$
$$|j-j'|\geq 2\Rightarrow\textrm{Supp}\,\ \varphi(2^{-j}\cdot)\cap \textrm{Supp}\,\ \varphi(2^{-j'}\cdot)=\emptyset,$$
$$j\geq 1\Rightarrow\textrm{Supp}\,\ \chi(\cdot)\cap \textrm{Supp}\,\ \varphi(2^{-j}\cdot)=\emptyset.$$
Then, we can define the nonhomogeneous dyadic blocks $\Delta_j$ and nonhomogeneous low frequency cut-off operator $S_j$ as follows:
$$\Delta_j{u}= 0,\,\ if\,\ j\leq -2,\quad
\Delta_{-1}{u}= \chi(D)u=\mathcal{F}^{-1}(\chi \mathcal{F}u),$$
$$\Delta_j{u}= \varphi(2^{-j}D)u=\mathcal{F}^{-1}(\varphi(2^{-j}\cdot)\mathcal{F}u),\,\ if \,\ j\geq 0,$$
$$S_j{u}= {\sum\limits_{j'=-\infty}^{j-1}}\Delta_{j'}{u}.$$

\begin{definition}[\cite{Bahouri 2011}]\label{de2.3}
Let $s\in\mathbb{R}$ and $1\leq p,r\leq\infty$. The nonhomogeneous Besov space $B^s_{p,r}(\R^d)$ consists of all tempered distribution $u$ such that
\begin{align*}
||u||_{B^s_{p,r}(\R^d)}\triangleq \Big|\Big|(2^{js}||\Delta_j{u}||_{L^p(\R^d)})_{j\in \Z}\Big|\Big|_{\ell^r(\Z)}<\infty.
\end{align*}
\end{definition}

Then, we have the following product and multiplier estimates, which will play an important role in the estimates of nonlocal term $\mathbf{P}(u)$ and convective term $u^Q\partial_xu$.

\begin{lemma}(\cite{Bahouri 2011})\label{lem2.1}
 (1) Algebraic properties: $\forall s>0,$ $B_{p, r}^s(\mathbb{R}^d)$ $\cap$ $L^\infty(\mathbb{R}^d)$ is a Banach algebra. $B_{p, r}^s(\mathbb{R}^d)$ is a Banach algebra $\Leftrightarrow B_{p, r}^s(\mathbb{R}^d)\hookrightarrow L^\infty(\mathbb{R}^d)\Leftrightarrow s>\frac{d}{p}$ or $s=\frac{d}{p},$ $r=1$.\\
 (2) For any $s>0$ and $1\leq p,r\leq\infty$, there exists a positive constant $C=C(d,s,p,r)$ such that
$$\|uv\|_{B^s_{p,r}(\mathbb{R}^d)}\leq C\Big(\|u\|_{L^{\infty}(\mathbb{R}^d)}\|v\|_{B^s_{p,r}(\mathbb{R}^d)}+\|v\|_{L^{\infty}(\mathbb{R}^d)}\|u\|_{B^s_{p,r}(\mathbb{R}^d)}\Big).$$
(3) Let $m\in \mathbb{R}$ and $f$ be an $S^m-$ multiplier (i.e., $f: \mathbb{R}^d\rightarrow \mathbb{R}$ is smooth and satisfies that $\forall \alpha\in \mathbb{N}^d$, there exists a constant $\mathcal{C}_\alpha$ such that $|\partial^\alpha f(\xi)|\leq \mathcal{C}_\alpha(1+|\xi|)^{m-|\alpha|}$ for all $\xi \in \mathbb{R}^d$). Then the operator $f(D)$ is continuous from $B_{p, r}^s$ to $B_{p, r}^{s-m}$.\\
(4) Let  $1\leq p, r\leq \infty$ and $s>\max\{1+\frac{1}{p}, \frac{3}{2}\}$. Then  we have
$$\|uv\|_{B_{p, r}^{s-2}(\mathbb{R})}\leq C\|u\|_{B_{p, r}^{s-2}(\mathbb{R})}\|v\|_{B_{p, r}^{s-1}(\mathbb{R})}.$$
Hence, for the terms $\mathbf{P}(u)$ and $\mathbf{P}(v)$, we have
\begin{eqnarray*}
  \|\mathbf{P}(u)-\mathbf{P}(v)\|_{B_{p, r}^{s-1}(\mathbb{R})}&\lesssim & \|u-v\|_{B_{p, r}^{s-1}(\mathbb{R})}\|u,v\|_{B_{p, r}^{s-1}(\mathbb{R})}^{Q-1}\|u, v\|_{B_{p, r}^s(\mathbb{R})},\\
 \|\mathbf{P}(u)-\mathbf{P}(v)\|_{B_{p, r}^s(\mathbb{R})}&\lesssim & \|u-v\|_{B_{p, r}^s(\mathbb{R})}\|u, v\|_{B_{p, r}^s(\mathbb{R})}^Q,\\
 \|\mathbf{P}(u)\|_{B_{p, r}^{s+1}(\mathbb{R})}&\lesssim & \|u\|_{B_{p, r}^{s+1}(\mathbb{R})}\|u\|_{B_{p, r}^s(\mathbb{R})}^Q.
\end{eqnarray*}
\end{lemma}

\begin{lemma}(\cite{Bahouri 2011, Li-Yin2})\label{lem2.2}
Let $1\leq p,r\leq \infty$ and
$
\sigma> - \min\{\frac{1}{p}, 1-\frac{1}{p}\}.
$
There exists a constant $C=C(p,r,\sigma)$ such that for any smooth solution to the following linear transport equation:
\begin{equation*}
\quad \partial_t f+v\pa_xf=g,\quad \; f|_{t=0} =f_0.
\end{equation*}
We have
\begin{align}\label{ES2}
\sup_{s\in [0,t]}\|f(s)\|_{B^{\sigma}_{p,r}(\mathbb{R})}\leq Ce^{CV_{p}(v,t)}\Big(\|f_0\|_{B^\sigma_{p,r}(\mathbb{R})}
+\int^t_0\|g(\tau)\|_{B^{s}_{p,r}(\mathbb{R})}\dd \tau\Big),
\end{align}
with
\begin{align*}
V_{p}(v,t)=
\begin{cases}
\int_0^t \|\nabla v(s)\|_{B^{\frac{1}{p}}_{p,\infty}(\mathbb{R})\cap L^\infty(\mathbb{R})}\dd s,&\ \ \mathrm{if} \; \sigma<1+\frac{1}{p},\\
\int_0^t \|\nabla v(s)\|_{B^{\sigma}_{p,r}(\mathbb{R})}\dd s,&\ \ \mathrm{if} \; \sigma=1+\frac{1}{p} \mbox{ and } r>1,\\
\int_0^t \|\nabla v(s)\|_{B^{\sigma-1}_{p,r}(\mathbb{R})}\dd s, &\ \ \mathrm{if} \;\sigma>1+\frac{1}{p}\ \mathrm{or}\ \{\sigma=1+\frac{1}{p} \mbox{ and } r=1\}.
\end{cases}
\end{align*}
\end{lemma}

\section{Non-uniform continuous dependence}
\setcounter{equation}{0}

In this section, we will give the proof of Theorem \ref{the1.1}. Before proceeding further, we need to introduce several important propositions to show that the solution $\mathbf{S}_t(u_0)$ can be approximated by $u_0-t(u_0)^Q\partial_xu_0+t\mathbf{P}(u_0)$ in a small time  near $t=0$.

Firstly, we establish the estimates of the difference between the solution $S_t(u_0)$ and  initial data $u_0$ in different Besov norms. That is

\begin{proposition}\label{pro 3.1}
Assume that $s>\max\{1+\frac{1}{p}, \frac{3}{2}\}$ with $1\leq p, r\leq \infty$ and $||u_0||_{B^s_{p,r}}\lesssim 1$. Then under the assumptions of Theorem \ref{the1.1}, we have
\bbal
&||\mathbf{S}_{t}(u_0)-u_0||_{B^{s-1}_{p,r}}\lesssim t||u_0||^Q_{B^{s-1}_{p,r}}||u_0||_{B^{s}_{p,r}},
\\&||\mathbf{S}_{t}(u_0)-u_0||_{B^{s}_{p,r}}\lesssim t\big(||u_0||^{Q+1}_{B^{s}_{p,r}}+||u_0||^Q_{B^{s-1}_{p,r}}||u_0||_{B^{s+1}_{p,r}}\big),
\\&||\mathbf{S}_{t}(u_0)-u_0||_{B^{s+1}_{p,r}}\lesssim t\big(||u_0||^Q_{B^{s}_{p,r}}||u_0||_{B^{s+1}_{p,r}}+||u_0||^Q_{B^{s-1}_{p,r}}||u_0||_{B^{s+2}_{p,r}}\big).
\end{align*}
\end{proposition}
{\bf Proof}\quad For simplicity, denote $u(t)=\mathbf{S}_t(u_0)$. Firstly, according to the local well-posedness result \cite{Mi 2013, Yan 2014, 1Yan 2014}, there exists a positive time $T=T(||u_0||_{B^s_{p,r}},s,p,r)$ such that the solution $u(t)$ belongs to $\mathcal{C}([0, T];  B_{p, r}^s)$. Moreover, by Lemmas \ref{lem2.1}-\ref{lem2.2}, for all $t\in[0,T]$ and $\gamma\geq s-1$,  there holds
\bal\label{u-estimate}
||u(t)||_{B^\gamma_{p,r}}\leq C||u_0||_{B^\gamma_{p,r}}.
\end{align}

Now we shall estimate the different Besov norms of the term $u(t)-u_0$, which can be bounded by $t$ multiplying the corresponding Besov norms of initial data $u_0$.

For $t\in[0,T]$, using (4) with $v=0$ and product estimates (2) in Lemma \ref{lem2.1}, we have from \eqref{u-estimate} that
\bbal
||u(t)-u_0||_{B^s_{p,r}}
&\leq \int^t_0||\pa_\tau u||_{B^s_{p,r}} \dd\tau
\\&\leq \int^t_0||\mathbf{P}(u)||_{B^s_{p,r}} \dd\tau+ \int^t_0||u^Q \pa_xu||_{B^s_{p,r}} \dd\tau
\\&\lesssim t\big(||u||^{Q+1}_{L_t^\infty(B^{s}_{p,r})}+||u||^Q_{L_t^\infty(L^\infty)}||u_x||_{L_t^\infty(B^{s}_{p,r})}\big)
\\&\lesssim t\big(||u||^{Q+1}_{L_t^\infty(B^{s}_{p,r})}+||u||^Q_{L_t^\infty(B^{s-1}_{p,r})}||u_x||_{L_t^\infty(B^{s}_{p,r})}\big)
\\&\lesssim t\big(||u_0||^{Q+1}_{B^{s}_{p,r}}+||u_0||^Q_{B^{s-1}_{p,r}}||u_0||_{B^{s+1}_{p,r}}\big),
\end{align*}
where we have used that $B_{p, r}^{s-1}$ is an Banach algebra with $s-1>\max\{\frac{1}{p}, \frac{1}{2}\}$ in the fourth inequality.

Following the same procedure of estimates as above, we have
\bbal
||u(t)-u_0||_{B^{s-1}_{p,r}}
&\leq \int^t_0||\pa_\tau u||_{B^{s-1}_{p,r}} \dd\tau
\\&\leq \int^t_0||\mathbf{P}(u)||_{B^{s-1}_{p,r}} \dd\tau+ \int^t_0||u^Q \pa_xu||_{B^{s-1}_{p,r}} \dd\tau
\\&\lesssim t||u||^Q_{L_t^\infty(B^{s-1}_{p,r})}||u||_{L_t^\infty(B^{s}_{p,r})}
\\&\lesssim t||u_0||^Q_{B^{s-1}_{p,r}}||u_0||_{B^{s}_{p,r}},
\end{align*}
and
\bbal
||u(t)-u_0||_{B^{s+1}_{p,r}}
&\leq \int^t_0||\pa_\tau u||_{B^{s+1}_{p,r}} \dd\tau
\\&\leq \int^t_0||\mathbf{P}(u)||_{B^{s+1}_{p,r}} \dd\tau+ \int^t_0||u^Q \pa_xu||_{B^{s+1}_{p,r}} \dd\tau
\\&\lesssim t\big(||u||^Q_{L_t^\infty(B^{s}_{p,r})}||u||_{L_t^\infty(B^{s+1}_{p,r})}
+||u||^Q_{L_t^\infty(B^{s-1}_{p,r})}||u||_{L_t^\infty(B^{s+2}_{p,r})}\big)
\\&\lesssim t\big(||u_0||^Q_{B^{s}_{p,r}}||u_0||_{B^{s+1}_{p,r}}+||u_0||^Q_{B^{s-1}_{p,r}}||u_0||_{B^{s+2}_{p,r}}\big),
\end{align*}
where we have used $ \|u^Q\|_{B_{p, r}^{s+1}} \lesssim \|u\|_{B_{p, r}^{s-1}}^{Q-1}\|u\|_{B_{p, r}^{s+1}},$ which can be proved by recurrence method in the third inequality. In fact, using product law (2) of Lemma \ref{lem2.1}, we have
\begin{eqnarray*}
   \|u^Q\|_{B_{p, r}^{s+1}}&\lesssim& \|u\|_{L^\infty} \|u^{Q-1}\|_{B_{p, r}^{s+1}}+\|u\|_{B_{p, r}^{s+1}}\|u^{Q-1}\|_{L^\infty}\\
    &\lesssim& \|u\|_{B_{p, r}^{s-1}} \|u^{Q-1}\|_{B_{p, r}^{s+1}}+\|u\|_{B_{p, r}^{s+1}}\|u\|_{B_{p, r}^{s-1}}^{Q-1}\\
      &\lesssim& \|u\|_{B_{p, r}^{s-1}}(\|u\|_{B_{p, r}^{s-1}} \|u^{Q-2}\|_{B_{p, r}^{s+1}}+\|u\|_{B_{p, r}^{s+1}}\|u\|_{B_{p, r}^{s-1}}^{Q-2})+\|u\|_{B_{p, r}^{s+1}}\|u\|_{B_{p, r}^{s-1}}^{Q-1}\\
       &\lesssim& \|u\|_{B_{p, r}^{s-1}}^2 \|u^{Q-2}\|_{B_{p, r}^{s+1}}+\|u\|_{B_{p, r}^{s+1}}\|u\|_{B_{p, r}^{s-1}}^{Q-1}\\
        &\vdots&\\
        &\lesssim& \|u\|_{B_{p, r}^{s-1}}^{Q-1}\|u\|_{B_{p, r}^{s+1}}.
        \end{eqnarray*}
Thus, we finish the proof of Proposition \ref{pro 3.1}.

Since Proposition \ref{pro 3.1} fails for critical index $(s, p, r)=(\frac32, 2, 1)$ due to
the lack of the estimate of solutions in $B_{2, 1}^{\frac{3}{3}}$, it needs special treatment. That is

\begin{proposition}\label{pro 3.2}
Assume that $(s,p,r)=(\frac32,2,1)$ and $||u_0||_{B^s_{p,r}}\lesssim 1$. Under the assumptions of Theorem \ref{the1.1}, we have
\bbal
&||\mathbf{S}_{t}(u_0)-u_0||_{L^\infty}\lesssim t||u_0||^{Q+1}_{C^{0,1}},
\\&||\mathbf{S}_{t}(u_0)-u_0||_{B^{s}_{p,r}}\lesssim t\big(||u_0||^{Q+1}_{B^{s}_{p,r}}+[E(u_0)]^Q||u_0||_{B^{s+1}_{p,r}}\big),
\\&||\mathbf{S}_{t}(u_0)-u_0||_{B^{s+1}_{p,r}}\lesssim t\big(||u_0||^Q_{B^{s}_{p,r}}||u_0||_{B^{s+1}_{p,r}}+[E(u_0)]^Q||u_0||_{B^{s+2}_{p,r}}\big),
\end{align*}
where $E(u_0)=||u_0||_{L^\infty}+||u_0||^{Q+1}_{C^{0,1}}$.
\end{proposition}

{\bf Proof}\quad For simplicity, we denote $u(t)=\mathbf{S}_t(u_0)$. Moreover, by Lemmas \ref{lem2.1}-\ref{lem2.2}, there holds for all $t\in[0,T]$ and $\gamma\geq \frac32$ that
\bal\label{u-estimate-2}
||u(t)||_{B^\gamma_{2,1}}\leq C||u_0||_{B^\gamma_{2,1}}, \qquad ||u(t)||_{C^{0,1}}\leq C||u_0||_{C^{0,1}}.
\end{align}

It follows by differential mean value theorem and the Minkowski inequality that
\bbal
||u(t)-u_0||_{L^\infty}
&\lesssim \int^t_0||\pa_\tau u||_{L^\infty} \dd\tau
\\&\lesssim \int^t_0||\mathbf{P}(u)||_{L^\infty} \dd\tau+ \int^t_0||u^Q \pa_xu||_{L^\infty} \dd\tau
\\&\lesssim t||u||^{Q+1}_{L_t^\infty(C^{0,1})}
\\&\lesssim t||u_0||^{Q+1}_{C^{0,1}},
\end{align*}
where we have used the fact that the operator $(1-\partial_x^2)^{-1}$ coincides with convolution by the function $x\rightarrow \frac{1}{2}e^{-|x|}$, and it can be proved by the Young inequality and H\"{o}lder inequality that
$$||\mathbf{P}(u)||_{L^\infty}\lesssim ||u||^{Q+1}_{L_t^\infty(C^{0,1})}.$$
Hence, we have
\bal\label{u-infinity}
||u(t)||_{L^\infty}\lesssim E(u_0).
\end{align}
Using the multiplier and product
estimates in Lemma \ref{lem2.1}, we obtain that
\bbal
||u(t)-u_0||_{B^{\frac32}_{2,1}}
&\leq \int^t_0||\pa_\tau u||_{B^{\frac32}_{2,1}} \dd\tau
\\&\leq \int^t_0||\mathbf{P}(u)||_{B^{\frac32}_{2,1}} \dd\tau+ \int^t_0||u^Q \pa_xu||_{B^{\frac32}_{2,1}} \dd\tau
\\&\leq Ct\big(||u||^{Q+1}_{L_t^\infty(B^{\frac32}_{2,1})}+||u||^Q_{L_t^\infty(L^\infty)}||u_x||_{L_t^\infty(B^{\frac32}_{2,1})}\big)
,
\end{align*}
which along with \eqref{u-estimate-2}-\eqref{u-infinity} yield the second estimate.

Similarly,
\bbal
||u(t)-u_0||_{B^{\frac52}_{2,1}}
&\leq \int^t_0||\pa_\tau u||_{B^{\frac52}_{2,1}} \dd\tau
\\&\leq \int^t_0||\mathbf{P}(u)||_{B^{\frac52}_{2,1}} \dd\tau+ \int^t_0||u^Q \pa_xu||_{B^{\frac52}_{2,1}} \dd\tau
\\&\leq Ct\big(||u||^Q_{L_t^\infty(B^{\frac32}_{2,1})}||u||_{L_t^\infty(B^{\frac52}_{2,1})}
+||u||^Q_{L_t^\infty(L^\infty)}||u||_{L_t^\infty(B^{\frac72}_{2,1})}\big)
\\&\leq Ct\big(||u_0||^Q_{B^{\frac32}_{2,1}}||u_0||_{B^{\frac52}_{2,1}}+[E(u_0)]^Q||u_0||_{B^{\frac72}_{2,1}}\big),
\end{align*}
which leads to the last inequality.
Thus, we complete the proof of Proposition \ref{pro 3.2}. $\Box$

With the different norms  estimates  of $u-u_0$ at hand, we have the following core estimates, which implies that for any bounded initial data $u_0$ in $B^s_{p,r}$, the corresponding solution $\mathbf{S}_t(u_0)$ can be approximated by $u_0-t(u_0)^Q\partial_xu_0+t\mathbf{P}(u_0)$ near $t=0$.
\begin{proposition}\label{pro 3.3}
Assume that $||u_0||_{B^s_{p,r}}\lesssim 1$. Then under the assumptions of Theorem \ref{the1.1}, there holds
\bbal
||\mathbf{S}_{t}(u_0)-u_0-t\mathbf{v}_0||_{B^{s}_{p,r}}\lesssim t^{2}\big(||u_0||^{Q+1}_{B^s_{p,r}}+||u_0||^Q_{B^{s-1}_{p,r}}||u_0||_{B^{s+1}_{p,r}}
+||u_0||^{2Q}_{B^{s-1}_{p,r}}||u_0||_{B^{s+2}_{p,r}}\big)
\end{align*} for $s>\max\{1+\frac{1}{p}, \frac{3}{2}\}$, $1\leq p, r\leq \infty,$
or
\bbal
||\mathbf{S}_{t}(u_0)-u_0-t\mathbf{v}_0||_{B^{s}_{p,r}}\lesssim t^{2}\big(||u_0||^{Q+1}_{B^s_{p,r}}+[E(u_0)]^Q||u_0||_{B^{s+1}_{p,r}}
+[E(u_0)]^{2Q}||u_0||_{B^{s+2}_{p,r}}\big)
\end{align*}
for $(s, p, r)=(\frac{3}{2}, 2, 1),$
here $\mathbf{v}_0=-u_0^Q\pa_x u_0+\mathbf{P}(u_0).$
\end{proposition}
{\bf Proof}\quad It should be noticed that according to (\ref{u-estimate})-(\ref{u-estimate-2}), $\|u\|_{B_{p, r}^s} \lesssim \|u_0\|_{B_{p, r}^s}\lesssim 1$ which will be frequently
used and be hidden later.

For $s>\max\{1+\frac{1}{p}, \frac{3}{2}\}$, $1\leq p, r\leq \infty,$ using (4) and product estimates (2) in Lemma \ref{lem2.1}, we obtain from  Propositions \ref{pro 3.1} that

\bbal
||u(t)-u_0-t\mathbf{v}_0||_{B^s_{p,r}}
&\leq \int^t_0||\pa_\tau u-\mathbf{v}_0||_{B^s_{p,r}} \dd\tau
\\&\leq \int^t_0||\mathbf{P}(u)-\mathbf{P}(u_0)||_{B^s_{p,r}} \dd\tau+\int^t_0||u^Q\pa_xu-u^Q_0\pa_xu_0||_{B^s_{p,r}} \dd\tau
\\&\leq \int^t_0||u(\tau)-u_0||_{B^s_{p,r}} \dd\tau+\int^t_0||u(\tau)-u_0||_{B^{s-1}_{p,r}} ||u(\tau)||_{B^{s+1}_{p,r}} \dd\tau
\\&\quad \ + \int^t_0||u(\tau)-u_0||_{B^{s+1}_{p,r}}  ||u_0||^Q_{B^{s-1}_{p,r}}\dd \tau
\\&\leq Ct^2\big(||u_0||^{Q+1}_{B^s_{p,r}}+||u_0||^Q_{B^{s-1}_{p,r}}||u_0||_{B^{s+1}_{p,r}}
+||u_0||^{2Q}_{B^{s-1}_{p,r}}||u_0||_{B^{s+2}_{p,r}}\big).
\end{align*}
Again using Lemma \ref{lem2.1} and Propositions \ref{pro 3.2}, for $(s, p, r)=(\frac{3}{2}, 2, 1),$ we arrive at
\bbal
||u(t)-u_0-t\mathbf{v}_0||_{B^{\frac32}_{2,1}}
&\leq \int^t_0||\pa_\tau u-\mathbf{v}_0||_{B^{\frac32}_{2,1}} \dd\tau
\\&\leq \int^t_0||\mathbf{P}(u)-\mathbf{P}(u_0)||_{B^{\frac32}_{2,1}} \dd\tau+\int^t_0||u^Q\pa_xu-u^Q_0\pa_xu_0||_{B^{\frac32}_{2,1}} \dd\tau
\\&\leq \int^t_0||u(\tau)-u_0||_{B^{\frac32}_{2,1}} \dd\tau+\int^t_0||u(\tau)-u_0||_{L^\infty} ||u_0,u(\tau)||^{Q-1}_{L^\infty}||u(\tau)||_{B^{\frac52}_{2,1}} \dd\tau
\\&\quad \ + \int^t_0||u(\tau)-u_0||_{B^{\frac52}_{2,1}}  ||u_0||^Q_{L^\infty}\dd \tau
\\&\lesssim t^2\big(||u_0||^{Q+1}_{B^s_{p,r}}+[E(u_0)]^Q||u_0||_{B^{s+1}_{p,r}}
+[E(u_0)]^{2Q}||u_0||_{B^{s+2}_{p,r}}\big).
\end{align*}

Thus, we complete the proof of Proposition \ref{pro 3.3}. $\Box$

Now, we move on the proof of Theorem \ref{the1.1}.

{\bf Proof of Theorem \ref{the1.1}}\quad  Let $\hat{\phi}\in \mathcal{C}^\infty_0(\mathbb{R})$ be an even, real-valued and non-negative funtion on $\R$ and satisfy
\begin{numcases}{\hat{\phi}(x)=}
1, &if $|x|\leq \frac{1}{4}$,\nonumber\\
0, &if $|x|\geq \frac{1}{2}$.\nonumber
\end{numcases}
Define the high frequency function $f_n$ and the low frequency function $g_n$ by
$$f_n=2^{-ns}\phi(x)\sin \bi(\frac{17}{12}2^nx\bi), \qquad g_n=\frac{12}{17}2^{-\frac nQ}\phi(x), \quad n\gg1.$$
It has been showed in \cite{Li 2020} that $\|f_n\|_{B_{p, r}^\sigma}\lesssim 2^{n(\sigma-s)}$.

Set $u^n_0=f_n+g_n$, consider Eq. (\ref{eq3}) with initial data $u_0^n$ and $f_n$, respectlively. Obviously, we have
\bbal
||u^n_0-f_n||_{B^s_{p,r}}=||g_n||_{B^s_{p,r}}\leq C2^{-\frac{n}{Q}},
\end{align*}
which means that
\bbal
\lim_{n\to\infty}||u^n_0-f_n||_{B^s_{p,r}}=0.
\end{align*}
For the case $s>\frac32$, it is easy to show that $||u^n_0,f_n||_{B^{s-1}_{p,r}}\lesssim 2^{-\frac nQ}$ and
\bbal
||u^n_0,f_n||_{B^{\sigma}_{p,r}}\leq C2^{(\sigma-s)n} \qquad  \mathrm{for} \qquad \sigma\geq s,
\end{align*}
which imply
\bbal
&\big(||u^n_0||^{Q+1}_{B^s_{p,r}}+||u^n_0||^Q_{B^{s-1}_{p,r}}||u^n_0||_{B^{s+1}_{p,r}}
+||u^n_0||^{2Q}_{B^{s-1}_{p,r}}||u^n_0||_{B^{s+2}_{p,r}}\big)\lesssim  1,
\\&\big(||f_n||^{Q+1}_{B^s_{p,r}}+||f_n||^Q_{B^{s-1}_{p,r}}||f_n||_{B^{s+1}_{p,r}}
+||f_n||^{2Q}_{B^{s-1}_{p,r}}||f_n||_{B^{s+2}_{p,r}}\big)\lesssim  1,
\end{align*}
For the case $s=\frac32$, it is easy to show that $E(u^n_0)+E(f_n)\lesssim 2^{-\frac nQ}$and
\bbal
||u^n_0,f_n||_{B^{\sigma}_{p,r}}\leq C2^{(\sigma-s)n} \qquad  \mathrm{for} \qquad \sigma\geq s,
\end{align*} 
which lead to
\bbal
&||u^n_0||^{Q+1}_{B^s_{p,r}}+[E(u^n_0)]^Q||u^n_0||_{B^{s+1}_{p,r}}
+[E(u^n_0)]^{2Q}||u^n_0||_{B^{s+2}_{p,r}}\lesssim  1,
\\&||f_n||^{Q+1}_{B^s_{p,r}}+[E(f_n)]^Q||f_n||_{B^{s+1}_{p,r}}
+[E(f_n)]^{2Q}||f_n||_{B^{s+2}_{p,r}}\lesssim  1.
\end{align*}

Furthermore, since $u_0^n$ and $f_n$ are both bounded in $B_{p, r}^s$, according to Proposition \ref{pro 3.3}, we deduce that
\bal\label{yyh}
\quad \ ||\mathbf{S}_{t}(u^n_0)-\mathbf{S}_{t}(f_n)||_{B^s_{p,r}}
\geq&~t\big|\big|(u^n_{0})^Q\pa_xu^n_{0}-(f_n)^Q\pa_xf_n-\mathbf{P}(u^n_0)+\mathbf{P}(f_n)\big|\big|_{B^s_{p,r}}-\big|\big|g_n\big|\big|_{B_{p, r}^s}-Ct^2\nonumber\\
\geq&~ t\big|\big|(u^n_{0})^Q\pa_xu^n_{0}-(f_n)^Q\pa_xf_n\big|\big|_{B^s_{p,\infty}}-C2^{-\frac nQ}-Ct^{2}.
\end{align}
Notice that
$$
(u^n_{0})^Q\pa_xu^n_{0}-(f_n)^Q\pa_xf_n=g^Q_n\pa_xf_n+(u^n_{0})^Q\pa_xg_n+\big((u^n_{0})^Q-f_n^Q-g^Q_n\big)\pa_xf_n,
$$
using Lemma \ref{lem2.1}, after simple calculation, we obtain
\bbal
||\big((u^n_{0})^Q-f_n^Q-g^Q_n\big)\pa_xf_n||_{B^s_{p,r}}&\leq ||\big((u^n_{0})^Q-f_n^Q-g^Q_n\big)||_{L^\infty}||f_n||_{B^{s+1}_{p,r}}\\
&\;\;+||\pa_xf_n||_{L^\infty}||\big((u^n_{0})^Q-f_n^Q-g^Q_n\big)||_{B^{s}_{p,r}}\\
&\leq C2^{-n(s-1)},\\
\big|\big|(u^n_{0})^Q\pa_xg_n\big|\big|_{B^s_{p,r}}&\leq ||u^n_0||^Q_{B^s_{p,r}}||g_n||_{B^{s+1}_{p,r}}\leq C2^{-\frac nQ},
\end{align*}
which together with fact
\begin{eqnarray*}
      \liminf_{n\rightarrow \infty} ||g^Q_n\pa_xf_n||_{B^s_{p,\infty}}\geq M
        \end{eqnarray*}
for some positive $M$ that has been showed in \cite{Li 2020} yield
\bbal
\liminf_{n\rightarrow \infty}||\mathbf{S}_t(f_n+g_n)-\mathbf{S}_t(f_n)||_{B^s_{p,r}}\gtrsim t\quad\text{for} \ t \ \text{small enough}.
\end{align*}
This completes the proof of Theorem \ref{the1.1}.

\vspace*{1em}
\noindent\textbf{Acknowledgements.}  J. Li is supported by the National Natural Science Foundation of China (Grant No.11801090) and Postdoctoral Science Foundation of China (Grant No.2020M672565). W. Zhu is partially supported by the National Natural Science Foundation of China (Grant No.11901092) and Natural Science Foundation of Guangdong Province (No.2017A030310634).

\end{document}